\documentclass[reqno, oneside, 12pt]{amsart}

\usepackage[letterpaper]{geometry}
\usepackage{enumerate, hyperref,url,amssymb,amsmath,amsthm,amsxtra,mathtools,mathrsfs,calc,nccmath,color}

\usepackage{calc}
\usepackage{graphicx}
\usepackage{caption}
\usepackage{subcaption}

\usepackage{needspace}

\usepackage{algorithm2e}

\SetStartEndCondition{ }{}{}%
\SetKwProg{Pydef}{def}{\string:}{}
\SetKwComment{tcc}{\#\ }{}
\SetKwFor{For}{for}{\string:}{}%
\SetKwFor{While}{while}{:}{fintq}%
\SetKwIF{If}{ElseIf}{Else}{if}{:}{elif}{else:}{}%
\SetKwFunction{Range}{range}
\SetKw{KwIn}{in}
\SetKw{KwNot}{not}
\SetKw{KwOr}{or}
\SetKw{KwInlineIf}{if}
\AlgoDontDisplayBlockMarkers\SetAlgoNoEnd\SetAlgoNoLine%

\newcommand{\Z}{\mathbb{Z}}
\newcommand{\Q}{\mathbb{Q}}

\newcommand{\ep}{\varepsilon}

\let\temp\phi
\let\phi\varphi
\let\varphi\temp
\newcommand{\texpdf}[2]{\texorpdfstring{#1}{#2}}

\newcommand{\SL}{\operatorname{SL}}
\newcommand{\GL}{\operatorname{GL}}

\newcommand{\sk}{\big|_k }

\newcommand{\pfrac}[2]{\left(\frac{#1}{#2}\right)}
\newcommand{\pmfrac}[2]{\left(\mfrac{#1}{#2}\right)}

\newcommand{\pMatrix}[4]{\left(\begin{matrix}#1 & #2 \\ #3 & #4\end{matrix}\right)}
\renewcommand{\pmatrix}[4]{\left(\begin{smallmatrix}#1 & #2 \\ #3 & #4\end{smallmatrix}\right)}

\newcommand{\tx}{\text}

\newtheorem{theorem}{Theorem}[section]

\newtheorem{proposition}[theorem]{Proposition}

\theoremstyle{remark}

\numberwithin{equation}{section}

\AtBeginDocument{%
   \def\MR#1{}
}

\begin{document}


\title{Scarcity of partition congruences on semiprime progressions}

\date{\today}
\author{Scott Ahlgren}
\address{Department of Mathematics\\
University of Illinois\\
Urbana, IL 61801} 
\email{sahlgren@illinois.edu} 

\author{Olivia Beckwith}
\address{Mathematics Department\\
Tulane University \\
New Orleans, LA 70087} 
\email{obeckwith@tulane.edu}

\thanks{The first author was  supported by a grant from the Simons Foundation (\#963004 to Scott Ahlgren).} \thanks{The second author was supported by NSF grant DMS-2401356 and by a grant from the Simons Foundation (\#953473 to Olivia Beckwith). 
}

 
\begin{abstract}  
In  recent work with Raum  the authors considered congruences for the ordinary partition function $p(n)$ of the  form $p(\ell Q^r n+\beta)\equiv 0\pmod\ell$ where $\ell, Q\geq 5$ are prime and  $r\in \{1,2\}$, and proved a number of results which show that such congruences are  scarce in a precise sense.
Here we  improve one of our results when $r=1$; in particular  we prove (outside of trivial cases) that  the set of primes $Q$ such that there exists $\beta\in \Z$ with $p(\ell Q n+\beta)\equiv 0\pmod \ell$ for all $n$ has density zero.  The proof involves a modification of part of our previous argument and an application of a recent theorem of Dicks regarding modular forms of half-integral weight and level one modulo $\ell$. 

 \end{abstract}

\dedicatory{Dedicated to Krishnaswami Alladi with our  thanks for his many years of tireless service to the mathematical community.}


\maketitle

 \section{Introduction}
 We consider linear congruences for the partition function, where $\ell$ is a prime number:
\begin{gather}\label{eq:generic-cong}  p(mn+\beta)\equiv 0\pmod\ell\quad\text{for all $n$}.
 \end{gather}
The theory of such congruences has a rich history which goes back to the three famous congruences of Ramanujan:
\begin{gather}\label{eq:ramcong}
p(\ell n+\beta)\equiv 0\pmod\ell
\quad\text{for }\ell=5, 7, 11
\tx{,}\quad \beta=\tfrac1{24}\pmod\ell.
\end{gather}
We briefly describe some of the important results which shape our understanding of such congruences.
Radu \cite{Radu_subbarao} proved a  conjecture of Subbarao by showing that there are no congruences of the form \eqref{eq:generic-cong} when $\ell = 2,3$. In later work \cite{Radu_aoconj}, Radu proved a conjecture of the first author and Ono by showing that for any congruence \eqref{eq:generic-cong}, we have
 $\ell\mid m$ and $\pfrac{1-24\beta}\ell\in \{0, -1\}$.

Let $\ell\geq 5$ be prime.  We will consider congruences of the form 
\begin{gather}\label{eq:partcong}
 p(\ell Q^r n+\beta)\equiv 0\pmod\ell, \quad r=0, 1, 2, \dots.
\end{gather}
 When $r=0$, work  of Kiming and Olsson \cite{Kiming-Olsson}
 and the first author and Boylan \cite{Ahlgren-Boylan}
 shows that the three congruences \eqref{eq:ramcong} are the only examples.
At the other extreme, for any $r\geq 4$ the work of Ono \cite{Ono_annals} and its generalizations 
\cite{Ahlgren_mathann, ahlgren-ono}
show that there are  infinitely many primes $Q$ which lead to \eqref{eq:partcong}.  The same is true generally for weakly holomorphic modular forms by 
the work of Treneer \cite{Treneer_1, Treneer_2}.

The situation is more interesting when $r=3$; here Atkin \cite{atkin_multiplicative} proved  a handful of 
examples for primes $\ell\leq 31$ (Johansson  \cite{johansson-2012}
added $22$ billion examples for these primes).  Recently the first author with Allen and Tang \cite{ahlgren-allen-tang} has showed that there are infinitely many such congruences for any prime $\ell\geq 5$  (this has been  generalized to a wide class of modular forms in \cite{AAD}).

It remains to consider congruences \eqref{eq:partcong} when $r=1, 2$. No examples are known in these cases, 
and in   
 work with Raum \cite{ahlgren-beckwith-raum}
 we proved a number of results which show that such congruences can exist only in very particular circumstances (if they exist at all). Here we briefly describe only some of those results.
Our basic approach
 involves analyzing  certain half-integral weight modular forms $f_{\ell, \delta}$
 (see \eqref{eq:fldeldef} below) with 
\begin{gather*}
f_{\ell, \delta} 
\equiv\sum_{\pfrac{-n}\ell=\delta}p\pmfrac{n+1}{24}q^\frac n{24}\pmod\ell,\quad \delta\in \{0,-1\}.
\end{gather*}

When  $r=2$, we showed \cite[Thm.~1.6]{ahlgren-beckwith-raum} 
that the congruence 
\eqref{eq:partcong} (outside of trivial cases) implies the restriction $Q \mid\mid 24 \beta -1 $, 
as well as  congruence relations involving 
$f_{\ell,\delta} \mid U_{Q^2}$, $f_{\ell,\delta} \mid V_{Q^2}$, and $f_{\ell,\delta} \otimes \chi_Q$ with $\delta = \left(\frac{1-24\beta}{\ell} \right)$. These relations enabled us to computationally rule out instances of \eqref{eq:partcong} with $r=2$ for  $17 \le \ell \le 1,000$ and $5 \le Q \le 10,000$.

When $r=1$  \cite[Thm. 1.4]{ahlgren-beckwith-raum} shows that either 
the set $S$ of primes $Q$ for which there is a congruence \eqref{eq:partcong} has density 0,
or the modular form $f_{\ell, \delta}$
(with 
$\delta = \left(\frac{1-24\beta}{\ell} \right)$) has very particular properties.  For example, 
$f_{\ell, \delta}$ is annihilated mod $\ell$ by all  Hecke operators $T_{Q^2}$ with $Q \equiv -1 \pmod{\ell}$.
This allows us to show that $S$ has density zero for $17\leq \ell<10,000$.
In \cite{ahlgren-beckwith-raum} we also show that 
 if there is a congruence \eqref{eq:partcong} 
with $r=1$ at $Q$ then the 
 Fourier expansion of $f_{\ell,\delta}$ 
 vanishes modulo $\ell$ on one of the square classes modulo $Q$.  This opens the door to ruling out $Q$ using an efficient  sieve-like computation. As a result we show that 
there is no such congruence (apart from the trivial examples arising from \eqref{eq:ramcong})
for $\ell<10,000$ and $Q<10^9$.

Our goal here is  to improve 
\cite[Thm. 1.4]{ahlgren-beckwith-raum} 
by replacing a technical part of the proof
with an alternate argument 
which allows us to rule out the second possibility described above.
This argument 
relies on a   recent result of Dicks \cite{robert-dicks} which
describes precisely  the situations in which a level one half-integral weight modular form can have its Fourier expansion supported on finitely many square classes modulo $\ell$.
We prove the following:
\begin{theorem} \label{thm:mainpart}
Suppose that $\ell\geq 5$ is prime,
and if  $\ell\in \{5, 7, 11\}$ suppose further that $\beta\not\equiv 1/24\pmod\ell$.
Then the set of primes $Q$ for which there exists a congruence of the form 
\begin{gather}\label{eq:Qellcong}
p(\ell Q n + \beta ) \equiv 0 \pmod{\ell}\quad \text{for all $n$}
\end{gather}
has density 0.
\end{theorem}

In Section 2 we give some brief  background material and we describe the result of Dicks as well as several fundamental results which limit the existence of a possible congruence \eqref{eq:partcong} with $r=1$.
 In Section 3, we use these results to give a proof of Theorem \ref{thm:mainpart}.

\section{Background}

\subsection{Modular forms with the eta multiplier}
The Dedekind eta function is defined by
\begin{gather*}
  \eta(\tau):=q^\frac1{24}\prod_{n=1}^\infty(1-q^n)
\tx{,}
\end{gather*}
where we use the notation 
\begin{gather*}
q:=e(\tau)=e^{2\pi i\tau},\  \ \ \tau \in \mathbb{H}.
\end{gather*}
The  eta function is a modular form of weight $1/2$ on $\SL_2(\Z)$; in particular, there is a multiplier $\nu_\eta$ with 
\begin{gather}\label{eq:etamult}
\eta(\gamma\tau)=\nu_\eta(\gamma)(c\tau+d)^\frac12\,\eta(\tau), \qquad \gamma=\pMatrix abcd\in \SL_2(\Z)\tx{.}
\end{gather}
where $\nu_{\eta}(\gamma)$ is the  24th root of unity given by the explicit formula in \cite[\S 4.1]{knopp}. 

If $f$ is a function on the upper half-plane,  $k\in \Z/2$,  and $\gamma=\pmatrix abcd\in \GL_2^+(\Q)$, 
we define
\begin{gather*}
\left(f\sk\gamma\right)(\tau):=(\det\gamma)^\frac k2(c\tau+d)^{-k}f(\tau).
\end{gather*}
Throughout the paper, $\ell\geq 5$ will denote a fixed prime number.
Given  $k\in \Z/2$,  a positive integer $N$, and a multiplier system $\nu$ on $\Gamma_0(N)$, we denote by $S_k (N,  \nu )$ the space of  cusp forms of weight $k$ and multiplier $\nu$ on $\Gamma_0(N)$ 
whose Fourier coefficients are  algebraic numbers which are integral at all primes above $\ell$.
Forms $f$ in this space  satisfy the transformation law
\begin{gather*}
f\sk \gamma=\nu (\gamma) f\ \ \ \text{for}\ \  \  \gamma=\pmatrix abcd\in \Gamma_0(N),
\end{gather*}
 vanish at the cusps of $\Gamma_0(N)$, 
 and have 
 Fourier expansions of the form
\begin{equation*}
f(\tau)=\sum_{n
    \equiv r\pmod{24}}a(n)q^\frac n{24}.
\end{equation*}
For $Q\geq 1$ we have the operators $U_Q$ and $V_Q$  whose action on $f=\sum a(n)q^\frac n{24}$ is given by 
\begin{gather*}
  f\big|U_Q:=\sum a(Qn)q^\frac n{24}
\quad\tx{and}\quad
  f\big|V_Q:=\sum a(n)q^\frac {Qn}{24}
\end{gather*}

\subsection{Fundamental results on congruences}\label{sec:prev_result}
As mentioned above, our method involves  analyzing  modular forms $f_{\ell, \delta}$ whose Fourier coefficients modulo $\ell$ capture  values of $p(n)$ in 
square-classes. To be precise, we showed in  (1.8)-(1.9) of \cite{ahlgren-beckwith-raum} that for $\delta\in \{0, -1\}$, there is a cusp form 
\begin{gather}\label{eq:fldelprop}
f_{\ell, \delta}\in 
\begin{cases}
S_{\frac{\ell^2 - 2\ell}{2}} \left( 1, \nu_\eta^{-1}\right)\ \ &\text{if $\delta=0$}\tx{,}\\
S_{\frac{\ell^2 - 2}{2}} \left( 1, \nu_\eta^{-1}\right)\ \ &\text{if $\delta=-1$}\tx{,}
\end{cases}
\end{gather}
with 
\begin{gather}\label{eq:fldeldef}
f_{\ell, \delta}=\sum a_{\ell,\delta}(n)q^\frac n{24} 
\equiv\sum_{\pfrac{-n}\ell=\delta}p\pmfrac{n+1}{24}q^\frac n{24}\pmod\ell.
\end{gather}

We  recall several important results. The first is  a theorem of Radu 
which shows that 
congruences for $p(n)$ propagate across square classes. Given a positive integer $m$ and an integer $\beta$, define the set 
\begin{gather*}
S_{m, \beta}:=\{\beta'\!\!\pmod m\ : \ 24\beta'-1\equiv a^2(24\beta-1)\!\!\pmod m\ \ \text{for some $a$ with $(a, 6m)=1$}\}.
\end{gather*}
Then we have
\begin{theorem}{\cite[Thm.~5.4]{Radu_subbarao}}\label{thm:raduclasses}
Suppose that $m$ and $\ell$ are positive integers and that for some integer $\beta$  we have a congruence $p(mn+\beta)\equiv 0\pmod\ell$.
Then for all $\beta'\in S_{m, \beta}$ we have a congruence $p(mn+\beta')\equiv 0\pmod\ell$.
\end{theorem}

The next  result of Radu shows that congruences can only occur on arithmetic progressions in certain square classes.
\begin{theorem}{\cite[Thm. 1.1]{Radu_aoconj}}\label{thm:radu_facts}
Suppose that $\ell\geq 5$ is prime.  Then a congruence \eqref{eq:Qellcong} can occur only if 
$\left(\frac{1-24\beta}\ell\right)\in\{0, -1\}$.
\end{theorem}
Finally, we require one of the main results of \cite{ahlgren-beckwith-raum}.
\begin{theorem}{\cite[Thm 1.5]{ahlgren-beckwith-raum}}\label{thm:abr_thm1.5}
Suppose that $\ell, Q\geq 5$ are distinct primes.
\begin{enumerate}
\item A congruence \eqref{eq:Qellcong} cannot occur with $\left(\frac{1-24\beta}Q\right)=0$.

\item  Fix $\delta\in \{0, -1\}$ and  $\ep\in \{\pm 1\}$.  If  there is a congruence of the form \eqref{eq:Qellcong}
with 
\begin{gather*}
 \pfrac{1-24\beta}\ell=\delta\ \ \text{and}\ \ \pfrac{1-24\beta}Q=\ep
\tx{,}
\end{gather*}
then
\begin{gather*}\label{eq:uqvq}
f_{\ell, \delta}\big| U_Q\equiv -\ep\pfrac{12}Q Q^{-1} \,f_{\ell, \delta}\big|V_Q\pmod\ell.
\end{gather*}

\end{enumerate}
\end{theorem}

\subsection{A mod \texpdf{$\ell$}{ell}
version of a theorem of Vign\'eras}
We recall a  recent theorem of Dicks which can be 
viewed as  a mod $\ell$, level one version of a theorem of Vign\'eras \cite{vigneras} (the result of Vign\'eras guarantees that any half integral weight modular form whose  Fourier expansion is supported on  finitely many square classes must be a linear combination of theta series). 
Dicks \cite{robert-dicks} gives a precise mod $\ell$ version of this result for half integral weight modular forms of level one whose weights  are sufficiently small with respect to $\ell$.  Here we state the result only in the case when the coefficient ring is $\Z$.
\begin{theorem}\label{thm:Robert} \cite[Thm. 1.3]{robert-dicks}
Suppose that $\ell \ge 5$ is prime 
and that $\lambda$ is a non-negative integer with $\lambda + \frac{1}{2} < \frac{\ell^2}{2}$. Suppose that $r$ is a positive integer with $(r,6) = 1$ and that $f \in S_{\lambda + \frac{1}{2}} (1, \nu_{\eta}^r) \cap \Z[\![q^{\frac{r}{24}} ]\!]$ satisfies
\begin{equation*}
    f \equiv \sum_{i=1}^m \sum_{n=1}^{\infty} a(t_i^2 n^2) q^{\frac{t_i n^2}{24}} \not\equiv 0 \pmod{\ell},
\end{equation*}
where each $t_i$ is a square-free positive integer. Then one of the following is true. 
\begin{enumerate}
\item $f \equiv a(1) \sum_{n=1}^{\infty} \left( \frac{12}{n} \right) n^{\lambda} q^{\frac{n^2}{24}} \pmod{\ell}$. In this case $r \equiv 1 \pmod{24}$ and $\lambda$ is even.
\item $f \equiv a(\ell) \sum_{n=1}^{\infty} \left( \frac{12}{n} \right) q^{\frac{\ell n^2}{24}} \pmod{\ell}.$ In this case $r \equiv \ell \pmod{24}$ and $\lambda \equiv \frac{\ell-1}{2} \pmod{\ell-~1}$. 
\item $f \equiv a(1)\sum_{n=1}^{\infty} \left( \frac{12}{n} \right) n^{\lambda} q^{\frac{n^2}{24}} +  a(\ell) \sum_{n=1}^{\infty} \left( \frac{12}{n} \right) q^{\frac{\ell n^2}{24}} \pmod{\ell}$, where $a(1) \not\equiv 0 \pmod{\ell}$ and $a(\ell) \not\equiv 0 \pmod{\ell}$. In this case, $r \equiv \ell \equiv 1 \pmod{24}$ and $\lambda \equiv \frac{\ell - 1}{2} \pmod{ \ell -1}$.
\end{enumerate}
\end{theorem}

\section{Proof of Theorem~\ref{thm:mainpart}}
By  the results of \cite{Kiming-Olsson} and \cite{Ahlgren-Boylan} we  have 
\begin{gather}
\label{eq:fl0}
f_{\ell, \delta}\equiv 0\pmod\ell\, \iff\, \ell\in \{5, 7, 11\}\quad \text{and} \quad \delta=0.
\end{gather}
From Theorem~\ref{thm:radu_facts} and the first part of Theorem~\ref{thm:abr_thm1.5} 
we see that in order to prove 
Theorem~\ref{thm:mainpart}  it  will suffice to establish the following:
\begin{proposition}\label{prop:reduce} Suppose that  $\ell\geq 5$ is prime and that 
$(\delta, \ep)$ is a pair with   $\delta\in \{0, -1\}$ and  $\ep\in \{\pm1\}$.
Then, outside of the three exceptional cases in \eqref{eq:fl0}, the set of primes $Q$ for which there is a congruence \eqref{eq:Qellcong} with  $\left(\frac{1-24\beta}\ell\right)=\delta$ and 
$\left(\frac{1-24\beta}Q\right)=\ep$ 
has density zero. 
\end{proposition}
To prove the proposition, fix 
such
a choice of  
$(\delta, \ep)$ 
 and let $f_{\ell, \delta}\not\equiv0\pmod\ell$ be as  in \eqref{eq:fldelprop}-- \eqref{eq:fldeldef}. Let $P_{\delta, \ep}$ be the set of primes $Q\ge 5$ for which there is a congruence \eqref{eq:Qellcong} for some $\beta$ with 
$\left( \frac{1-24 \beta }{\ell} \right)  = \delta$ and $\left( \frac{1-24 \beta }{Q} \right)  = \ep$. Our goal is to show that $P_{\delta, \ep}$ has density 0, i.e., that
\begin{equation*}
    \lim_{X \to \infty} \frac{ \# \{ Q \in P_{\delta, \ep} : Q < X \}}
{X/\log X} = 0.
\end{equation*}

Let $Q \in P_{\delta, \ep}$. From Theorem \ref{thm:raduclasses}, we have $p(\ell Q n + \beta') \equiv 0 \pmod{\ell}$ for all $\beta' \in S_{\ell Q, \beta}$ and $n \in \mathbb{Z}$.
Recalling the definition of $S_{\ell Q, \beta}$, we conclude that 
\begin{gather*}
  p(N)\equiv 0\pmod\ell\quad\text{if}\quad \pfrac{1-24N}\ell=\delta
  \quad\text{and}
  \quad   \pfrac{1-24N}Q=\ep,
\end{gather*}
or in other words that
\begin{gather*}
  p\pfrac{n+1}{24}\equiv 0\pmod\ell\quad\text{if}\quad \pfrac{-n}\ell=\delta
\quad\text{and}\quad   \pfrac{-n}Q=\ep.
\end{gather*}
From this we conclude that 
\begin{equation}\label{eq:FourierModEll}
f_{\ell, \delta} \equiv \sum_{\left( \frac{-n}{Q} \right) = -\ep} a_{\ell,\delta} (n) q^\frac n{24} + \sum a_{\ell, \delta} ( Qn) q^{\frac{Qn}{24}} \pmod{\ell}.
\end{equation}
On the other hand, from the second part of Theorem~\ref{thm:abr_thm1.5} we have
\begin{gather}\label{eq:uqvqcong}
    \sum a_{\ell, \delta} (Q n ) q^{\frac{n}{24}} \equiv -\ep\pfrac{12}Q Q^{-1} \sum a_{\ell, \delta}(n) q^{\frac{Q n}{24}} \pmod{\ell}.
\end{gather}
By \eqref{eq:uqvqcong} we see that 
 $a_{\ell, \delta} (Qn )\not\equiv0\pmod\ell\implies$  $Q \mid n$,
 which leads to
\begin{equation}\label{eq:FourierQMults}
 \sum a_{\ell, \delta} ( Qn) q^{\frac{Qn}{24}} \equiv  \sum a_{\ell, \delta} ( Q^2n) q^{\frac{Q^2n}{24}} \pmod{\ell}.
\end{equation}
Finally, combining \eqref{eq:FourierModEll} and \eqref{eq:FourierQMults}, we deduce that for every $Q \in P_{\delta, \ep}$, we have
\begin{equation*}
f_{\ell,\delta} \equiv \sum_{\left( \frac{-n}{Q} \right) = - \ep} a_{\ell, \delta}(n) q^{n/24} + \sum a_{\ell, \delta} (Q^2 n ) q^{Q^2 n/24}. 
\end{equation*}
In other words,
\begin{equation}\label{eq:an_nonvanish}
  a_{\ell,\delta}(n) \not\equiv 0 \pmod\ell \   \implies\ 
  \left( \frac{-n}{Q} \right) = -\ep \quad\text{ or }\quad Q^2 \mid n
\end{equation}

 We now claim that there is no 
  finite sequence of squarefree integers $t_i$ such that
 \begin{equation}
\label{eq:fin_many_classes}
f_{\ell, \delta} \equiv \sum_{i=1}^m \sum_{n =1}^\infty a_{\ell, \delta}(t_i n^2) q^{\frac{t_i n^2}{24}} \pmod{\ell}.
\end{equation}
To prove this we recall  that the weight $\lambda+1/2$ of $f_{\ell, \delta}$ is given in \eqref{eq:fldelprop} and is less than $\ell^2/2$, so  we may apply Theorem~\ref{thm:Robert}  to $f_{\ell, \delta}$. By \eqref{eq:fldelprop} we have $r = -1$, so cases (1) and (3) of Theorem~\ref{thm:Robert} are not possible. 
We also find that case (2) is not possible, since 
\begin{equation*}
\lambda=
\begin{cases}
\frac{\ell^2 - 2\ell - 1}2\quad&\text{if}\quad \delta=0,\\
\frac{\ell^2 - 3}2\quad&\text{if}\quad \delta=-1,
\end{cases}
\end{equation*}
and in both cases we have 
$\lambda  \equiv -1  \pmod{\ell -1}$.
 This establishes the claim.

It follows  that there is an infinite sequence of positive square free numbers $t_1<t_2<\cdots$ such that for each $i$ there exists some $n_i$ with $a_{\ell, \delta}(t_in_i^2)\not\equiv 0\pmod\ell$.
Suppose that $Q\in P_{\delta, \ep}$ and that $Q\mid n_i$. 
From the second part of Theorem~\ref{thm:abr_thm1.5} we have
\begin{gather*}
f_{\ell, \delta}\big| U_{Q^2}\equiv -\ep\pfrac{12}Q Q^{-1} \,f_{\ell, \delta}\pmod\ell.
\end{gather*}
So we can conclude that  $a_{\ell, \delta}\left(t_in_i^2/Q^2\right)\not\equiv 0\pmod\ell$.  Removing additional factors of $Q^2$ if necessary  we can assume that for all $Q\in P_{\delta, \ep}$ and all $i$ we have $Q\nmid n_i$.
This gives an infinite sequence $t_1<t_2<\cdots$ of positive squarefree integers such that 
\begin{equation}\label{eq:value_ep}
\pfrac{-t_i}Q=-\ep\quad\text{for all}\quad Q\in P_{\delta, \ep}.
\end{equation}
We may then choose a subsequence $t_{i_1}<t_{i_2}<\cdots$ with the property that
for every $k\geq 1$, $t_{i_k}$ is divisible by a prime $p_k$ which does not divide any of $t_{i_1},\dots, t_{i_{k-1}}$.

Fix $k\geq 1$ and consider the 
sum
\begin{equation}\label{eq:skdef}
    S_k(X):=\sum_{\substack{Q\ \tx{prime}\\ Q\leq X}}
\left(1-\ep\pfrac{-t_{i_1}}Q\right)\cdots\left(1-\ep\pfrac{-t_{i_k}}Q\right).
\end{equation}
From \eqref{eq:value_ep} we see that 
\begin{equation}\label{eq:sk_lowerbound}
    S_k(X)\geq 2^k\cdot  \#\left\{Q\in P_{\delta, \ep}\ : \ Q\leq X\right\}.
\end{equation}
 If $t$ is not a square then $\pfrac t\bullet$ corresponds to  a non-trivial Dirichlet character and  it follows from Dirichlet's theorem on primes in  progression that
\begin{equation*}
\sum_{\substack{Q\ \tx{prime}\\ Q\leq X}}\pfrac{t}Q=o\left(\frac X{\log X}\right).
\end{equation*}
Expanding out the products in \eqref{eq:skdef} and using this fact along with our assumption on the $t_{i_k}$ gives
\begin{equation*}
S_k(X)\sim  \frac{X}{\log X}.
\end{equation*}
Combining this with \eqref{eq:sk_lowerbound}
we obtain
\begin{equation*}
    \lim_{X\to\infty}\frac{\#\left\{Q\in P_{\delta, \ep}\ : \ Q\leq X\right\}}{X/\log X}\leq 2^{-k}.
\end{equation*}
Since $k$ is arbitrary this shows that $P_{\delta, \ep}$ has density zero. This establishes Proposition~\ref{prop:reduce} and with it Theorem~\ref{thm:mainpart}.

\bibliographystyle{amsalpha}
\bibliography{part_cong}

\end{document}